\documentclass[12pt]{iopart}
\usepackage{graphicx,latexsym,color}

\usepackage{cite}   
\usepackage[colorlinks]{hyperref}
\usepackage{xcolor}
\hypersetup{
    linkcolor={red!50!black},
citecolor={blue!50!black},
urlcolor={blue!20!black},
}

\usepackage[nolist]{acronym}
\usepackage[colorinlistoftodos,prependcaption,textsize=small]{todonotes}
\usepackage[normalem]{ulem}
\usepackage{regexpatch}
\makeatletter
\xpatchcmd{\@todo}{\setkeys{todonotes}{#1}}{\setkeys{todonotes}{inline,#1}}{}{}
\makeatother

\colorlet{colorMG}{blue}
\colorlet{colorBN}{green}
\colorlet{colorJT}{yellow}
\colorlet{colorLJ}{orange}
\colorlet{colorSN}{red}
\colorlet{colorYI}{teal!75}
\colorlet{colorAL}{blue!50}

\usepackage[T1]{fontenc}
\usepackage[swedish,english]{babel}
\usepackage[applemac]{inputenc} 

\usepackage{amsgen,amsfonts,amsbsy,amssymb}

\usepackage{bm}
\usepackage{tikz}
\usepackage{pgfplots}

\newcommand{\R}{\mathbb{R}}
\newcommand{\C}{\mathbb{C}}

\newcommand{\ie}{\textit{i.e.}\/, }
\newcommand{\eg}{\textit{e.g.}\/, }
\newcommand{\cf}{\textit{cf.}\/, }

\providecommand*{\mrm}[1]{\mathrm{#1}}
\providecommand*{\unit}[1]{\ensuremath{\mrm{\,#1}}}

\providecommand*{\iu}{\ensuremath{\mrm{i}}}

\newcommand{\toh}{\hat{\to}}

\newtheorem{theorem}{\em Theorem}[section]

\newtheorem{definition}{\em Definition}[section]
\newtheorem{example}{\em Example}[section]

\newcommand{\minimize}{\mrm{minimize}}

\newcommand{\subto}{\mrm{subject\ to}}

\renewcommand{\Re}{\ensuremath{\mrm{Re}}}	
\renewcommand{\Im}{\ensuremath{\mrm{Im}}}	

\def\Xint#1{\mathchoice
   {\XXint\displaystyle\textstyle{#1}}%
   {\XXint\textstyle\scriptstyle{#1}}%
   {\XXint\scriptstyle\scriptscriptstyle{#1}}%
   {\XXint\scriptscriptstyle\scriptscriptstyle{#1}}%
   \!\int}
\def\XXint#1#2#3{{\setbox0=\hbox{$#1{#2#3}{\int}$}
     \vcenter{\hbox{$#2#3$}}\kern-.5\wd0}}

\def\dashint{\Xint-}

\begin{document}

\title{Passive approximation and optimization using B-splines}

\author{Y.~Ivanenko$^1$, M.~Gustafsson$^2$,  B.L.G.~Jonsson$^3$, A.~Luger$^4$, B.~Nilsson$^5$, S.~Nordebo$^1$, J.~Toft$^5$}

\address{$^1$ Department of Physics and Electrical Engineering, Linn\ae us University, 351 95 V\"{a}xj\"{o}, Sweden.
E-mail: \{yevhen.ivanenko,sven.nordebo\}@lnu.se.}
\address{$^2$ Department of Electrical and Information Technology, Lund University, Box 118, 
221 00 Lund, Sweden. E-mail: mats.gustafsson@eit.lth.se.}
\address{$^3$ School of Electrical Engineering, KTH Royal Institute of Technology, 
100 44 Stockholm, Sweden. E-mail: ljonsson@kth.se.}
\address{$^4$ Department of Mathematics, Stockholm University
106 91 Stockholm, Sweden. E-mail: luger@math.su.se.}
\address{$^5$ Department of Mathematics, Linn\ae us University, 351 95 V\"{a}xj\"{o}, Sweden.
E-mail: \{borje.nilsson,joachim.toft\}@lnu.se.}

\maketitle

\begin{abstract}
A passive approximation problem is formulated where the target function is an arbitrary 
complex valued continuous function defined on an approximation domain consisting of a finite union of closed and bounded intervals on the real axis. 
The norm used is a weighted $\mrm{L}^p$-norm where $1\leq p\leq\infty$. 
The approximating functions are Herglotz functions generated by a measure 
with H\"{o}lder continuous density in an arbitrary neighborhood of the approximation domain. 
Hence, the imaginary and the real parts of the approximating functions are H\"{o}lder continuous functions given by the density of the measure and its Hilbert transform, respectively.
In practice, it is useful to employ finite B-spline expansions to represent the generating measure.
The corresponding approximation problem can then be posed as a finite-dimensional convex optimization problem which is amenable for numerical solution.
A constructive proof is given here showing that the convex cone of approximating functions generated by finite uniform B-spline expansions of fixed arbitrary order 
(linear, quadratic, cubic, etc) is dense in the convex cone of Herglotz functions which are locally H\"{o}lder continuous in a neighborhood of the approximation domain, as mentioned above.
As an illustration, a typical physical application example is included regarding the passive 
approximation and optimization of a linear system having metamaterial characteristics.  
\end{abstract}

\section{Introduction}


The basic assumption of passivity assigned to a time-translational invariant linear system implies many important physical constraints,
one of the most fundamental property being causality \cite{Zemanian1965}. It is well known that causality by itself limits the response of a linear system
via the associated Kramers-Kronig relations, see \eg \cite{Nussenzveig1972,King2009,Milton+Eyre1997,Haakestad+Skaar2005}.
However, it is important to note that the additional assumption of passivity may imply severe bandwidth limitations which are
not present if the system is merely assumed to be causal. 
Bandwidth limitations on passive systems can be derived based on the theory of Herglotz functions 
where a complexified angular frequency is the complex variable associated with the holomorphic Fourier transform \cite{Rudin1987}.
In particular, sum rules can be used to derive physical bounds on passive systems that can be useful in a 
variety of applications, see \eg \cite{Bernland+etal2011b}.
A classical example is with the bounds that were derived by Fano \cite{Fano1950} on broadband matching based on lossless networks.
More recently, sum rules have been used to derive physical bounds in applications such as with radar absorbers \cite{Rozanov2000}, 
high-impedance surfaces \cite{Gustafsson+Sjoberg2011}, passive metamaterials \cite{Gustafsson+Sjoberg2010a}, broadband quasi-static cloaking \cite{Cassier+Milton2017},
scattering~\cite{Sohl+etal2007a,Bernland+etal2011b}, antennas~\cite{Gustafsson+etal2010a,Jonsson+etal2013}, 
reflection coefficients~\cite{Gustafsson2010b}, waveguides~\cite{Vakili+etal2014}, and periodic structures~\cite{Gustafsson+etal2012c}, etc. 
Several other limitations have also been derived based on general  Herglotz function theory, such as
with the speed-of-light limitations in passive linear media described in \cite{Welters+etal2014}.

The application of sum rules to obtain bandwidth limitations on a physical system relies on the existence of
certain moments \cite{Akhiezer1965} as well as some a priori knowledge, or even measurements \cite{Vakili+etal2014},
regarding the low- and/or high-frequency asymptotic expansion coefficients of the corresponding Herglotz function \cite{Bernland+etal2011b}.
These parameters are typically related to the static and/or the high frequency responses of a material, or a structure. 
However, there are also many circumstances where sum rules do not apply, typically in situations when there are losses involved
that inhibit the necessary asymptotic expansions. A simple example is with a metal backed radar absorber when the metal is not perfectly conducting, see \eg \cite{Nordebo+etal2014b}.
Another more general example is with the passive approximation of an arbitrary continuous function on a closed interval on the real axis,
corresponding \eg to the approximate realization of some ``idealized'' linear system specified to have
some ``desired''  response over a given fixed frequency band.  
In a very special case, when the target function is the negative of some Herglotz function with a given high-frequency asymptotics, a
sum rule can be used to derive a lower bound on the approximation error in the max-norm, \cf \cite{Gustafsson+Sjoberg2010a,Ivanenko+etal2017}.
A typical physical application example for this problem formulation is with the passive approximation of a linear system having metamaterial characteristics.
In particular,  when a dielectric medium is specified to have a constant negative value of permittivity (an inductive property) over a given bandwidth, 
the passivity of the material will imply severe bandwidth limitations and the corresponding sum rule yields a 
lower bound on the max-norm of the approximation error, \cf \cite{Gustafsson+Sjoberg2010a}.
However, the lower bound on the permittivity error is not tight since it corresponds to a weighted norm on the corresponding Herglotz function.
Hence, in a general situation the sum rules may not give bounds that are tight for the specific problem at hand,
or there are simply no known sum rules that apply.  
An alternative is then to employ a passive approximation formulation based on convex optimization to study the 
 corresponding physical bounds \cite{Nordebo+etal2014b}. 

The purpose of the present paper is to provide a rigorous mathematical framework for the convex optimization approach based on general weighted $\mrm{L}^p$-norms.
In particular, it is proved that the set of generating measures given by finite uniform B-spline expansions of fixed arbitrary order yield
a set of approximating Herglotz functions which is dense in the set of Herglotz functions which are H\"{o}lder continuous in some  neighborhood of the approximation domain. The precise meaning of this statement is elaborated on throughout the paper
based on the theory of passive systems  \cite{Zemanian1965}, 
Herglotz functions (also known as Nevanlinna, Herglotz-Nevanlinna, Pick, R-, and positive real (PR) functions) \cite{Zemanian1965,Kac+Krein1974,Akhiezer1965,Nussenzveig1972,Gesztesy+Tsekanovskii2000}, 
B-splines \cite{Boor1968,Boor1972} and H\"{o}lder continuity \cite{Kress1999,King2009}.
Further details and application examples regarding passive approximation and optimization
using B-splines are also given in \cite{Nordebo+etal2014b,Nordebo+etal2017a,Ivanenko+Nordebo2016a,Ivanenko+etal2017}.

\section{General properties of Herglotz functions}

\subsection{The generating measure}
In the following, a complex number $z\in\C$ is written $z=x+\iu y$ with $x,y\in\R$.
A Herglotz function $h(z)$ is an analytic function defined on the open upper half-plane $\C^+=\{z\in\C|\Im\{z\}>0\}$ with the property that
$\Im\{h(z)\}\geq 0$ for $z\in\C^+$. 
It can be shown that $h(z)$ is a Herglotz function if and only if it can be represented as
\begin{equation}\label{eq:Herglotz1}
h(z)=b_1z+c+\int_{-\infty}^{\infty}\frac{1+\xi z}{\xi-z}{\rm d}\mu(\xi),
\end{equation}
where $b_1\geq 0$, $c\in \R$ and $\mu$ is a finite positive Borel measure, 
see \eg \cite{Beltrami+Wohlers1966,Kac+Krein1974,Akhiezer1965,Nussenzveig1972,Gesztesy+Tsekanovskii2000}.
It is also useful to introduce the positive Borel measure $\beta$ with ${\rm d}\beta(\xi)=(1+\xi^2){\rm d}\mu(\xi)$ so that
\begin{equation}\label{eq:Herglotz2}
h(z)=b_1z+c+\int_{-\infty}^{\infty}\left(\frac{1}{\xi-z}-\frac{\xi}{1+{\xi}^2}\right){\rm d}\beta(\xi),
\end{equation}
where $\int_{\R}{\rm d}\beta(\xi)/(1+\xi^2)<\infty$.  
In \eref{eq:Herglotz1} and \eref{eq:Herglotz2}, 
the constant $b_1$ is given by $b_1=\lim_{z\toh\infty}h(z)/z$ where $z\toh \infty$ 
means that $|z|\to \infty$ in the Stolz cone $\varphi\leq\arg z \leq\pi-\varphi$ for any $\varphi\in(0,\pi/2]$.
The constant $c$ is given by $c=\Re\left\{h(\iu)\right\}$. 
The positive measure $\beta$ is uniquely determined by the Herglotz function $h(z)$
from the Stieltjes inversion formula
\begin{equation}\label{eq:measdef}
\beta\left((x_1,x_2)\right)+\frac{1}{2}\beta\left(\{x_1\}\right)+\frac{1}{2}\beta\left(\{x_2\}\right)=\displaystyle\lim_{y\rightarrow 0+}\frac{1}{\pi}\int_{x_1}^{x_2}\Im\{h(\xi+\iu y)\}{\rm d}\xi, 
\end{equation}
including the possibility of having point masses at any $x_i\in\R$, see \cite{Kac+Krein1974,Gesztesy+Tsekanovskii2000}.

The following theorem further characterizes the connection between the generating measure and the imaginary part of the Herglotz function,
and it is stated here without proof, see \eg \cite[p.~7]{Kac+Krein1974}.
\begin{theorem}\label{theo:limitsofbetaandimh}  \rm
The following two limits are equivalent in the sense that the existence of one limit implies the existence of the other,
\begin{eqnarray}\label{eq:betaprimelimit}
f(x)=\displaystyle \lim_{\epsilon\rightarrow 0+}\frac{\beta((x-\epsilon,x+\epsilon))}{2\epsilon}, \vspace{0.2cm} \label{eq:betaprimelimitdef} \\
g(x)=\displaystyle \lim_{y\rightarrow 0+}\frac{1}{\pi}\Im\{h(x+\iu y)\}, \label{eq:glimitdef}
\end{eqnarray}
and $f(x)=g(x)$. Furthermore, if the limit $g(x)$ exists as a bounded function at all points $x$ of an interval $(x_1,x_2)$, then
the measure $\beta$ is absolutely continuous on that interval and the density $\beta^\prime$ (Radon-Nikod\'{y}m derivative) is given by $\beta^\prime(x)=g(x)$ almost everywhere.
\hfill $\blacksquare$
\end{theorem}
In the case when the measure is absolutely continuous it is common to adopt the notation ${\rm d}\beta(x)=\beta^\prime(x){\rm d}x=\frac{1}{\pi}\Im\{h(x+\iu 0)\}{\rm d}x$,
see \eg \cite{Nussenzveig1972,Gesztesy+Tsekanovskii2000}.

Symmetric Herglotz functions satisfy the symmetry requirement $h(z)=-h(-z^*)^*$, where $z\in\C^+$.
In this case, $\beta$ is an even measure and \eref{eq:Herglotz2} can be simplified as
\begin{equation}\label{eq:Herglotz3}
h(z)=b_1z+\int_{-\infty}^{\infty}\frac{1}{\xi-z}{\rm d}\beta(\xi),
\end{equation}
for $z\in\C^+$, and the integral is taken as a symmetric limit at infinity.

\subsection{Asymptotic properties and sum rules}\label{sect:asymptotics}
Herglotz functions have the general asymptotic behavior $h(z)=a_{-1}z^{-1}+o(z^{-1})$ as $z\toh 0$ where $a_{-1} \leq 0$,  
and $h(z)=b_1z+o(z)$ as $z\toh\infty$ where $b_{1} \geq 0$, see \eg \cite{Akhiezer1965,Bernland+etal2011b}.
Suppose now that the small and large argument asymptotic expansions of a symmetric Herglotz function $h(z)$ are given by 
\begin{equation}\label{eq:hLFHFassymptotic}
\hspace{-1cm} h(z)=\left\{\begin{array}{l}
a_{-1}z^{-1}+a_1z+\ldots+ a_{2N_0-1}z^{2N_0-1}+o(z^{2N_0-1})
 \quad \textrm{as}\ z\toh 0, \vspace{0.2cm}\\
b_1z+b_{-1}z^{-1}+\ldots+b_{1-2N_{\infty}}z^{1-2N_{\infty}}+ o(z^{1-2N_{\infty}})
 \quad \textrm{as}\ z\toh \infty, 
\end{array}\right.
\end{equation}
where all the expansion coefficients are 
real-valued, $N_0$ and $N_{\infty}$ are non-negative integers (or possibly infinity) and where $1-N_{\infty}\leq N_0$. 
In this case, it is possible to show that the following integral identities (sum rules) hold
\begin{equation}\label{eq:Herglotzidentity}
\hspace{-1cm} \frac{2}{\pi} \int_{0+}^{\infty}\frac{\Im\{h(x)\}}{x^{2k}}{\rm d} x\stackrel{\rm def}{=}
\lim_{\varepsilon\rightarrow 0+}\lim_{y\rightarrow 0+}\frac{2}{\pi}\int_{\varepsilon}^{1/\varepsilon}\frac{\Im\{h(x+\iu y )\}}{x^{2k}}{\rm d}x=a_{2k-1}-b_{2k-1},
\end{equation}
for $k = 1-N_{\infty},\ldots, N_0$, see \eg \cite{Akhiezer1965,Bernland+etal2011b}.

\subsection{Boundary values on the real axis}\label{sect:Holderspaces}
The boundary values of Herglotz functions on the real axis can generally be interpreted in the sense of 
tempered distributions and where the associated Hilbert transform is defined as a Cauchy principal value integral, see \eg  \cite{Zemanian1965,Beltrami+Wohlers1966,Nussenzveig1972,King2009}.
In the present application concerning passive approximation of functions that are
known to be continuous on finite intervals on the real line, it is natural to consider approximating functions with similar properties.
However, as the following example shows, it is not sufficient to assume
that the measure $\beta$ is absolutely continuous with continuous density $\beta^\prime$ on a given interval $\Omega$ to guarantee the continuity (or boundedness) of the Hilbert transform on $\Omega$.

\begin{example}\label{ex:nonholdercontinuous}\rm
Consider the function 
\begin{equation}\label{eq:exnonholdercontinuous}
\beta^\prime(x)=\left\{
\begin{array}{ll}
\displaystyle -\frac{1}{\ln x} & x\in (0,a], \\
0 & x\in [-a,0],
\end{array}\right.
\end{equation}
which is a continuous function on $\Omega=[-a,a]$ where $0<a<1$. 
It is readily seen that the Hilbert transform as a Cauchy principal value integral does not exist at $x=0$, see also \cite[ p.~97]{Kress1999},
and hence that the corresponding Herglotz function does not admit a continuous extension to $\Omega$.
\hfill $\square$
\end{example}

H\"{o}lder continuity is the adequate property to avoid the difficulties implied by the Example \ref{ex:nonholdercontinuous}.
A real or complex valued function $f(x)$ defined on the closed and bounded interval $\Omega$ is 
H\"{o}lder continuous (sometimes also referred to as uniformly H\"{o}lder continuous) with
H\"{o}lder exponent $0<\alpha< 1$ if there exists a constant $C$ such that
\begin{equation}\label{eq:Holdercontinuity}
\left| f(x)-f(y) \right|\leq C \left| x-y \right|^\alpha,
\end{equation}
for all $x,y\in\Omega$. The corresponding linear space $C^{0,\alpha}(\Omega)$ of H\"{o}lder continuous functions is a Banach space
with norm
\begin{equation}\label{eq:Holdernorm}
\|f\|_{\alpha}=\max_{x\in\Omega}|f(x)|+\sup_{x,y\in\Omega,x\neq y}\frac{\left| f(x)-f(y) \right|}{\left| x-y \right|^\alpha},
\end{equation}
see \cite[pp.~94-104]{Kress1999}. A smaller H\"{o}lder exponent implies a larger space, so that $C^{0,\alpha}(\Omega)\subset C^{0,\alpha^\prime}(\Omega)$,
where $\alpha^\prime<\alpha$.
Note that the example function in \eref{eq:exnonholdercontinuous} is not H\"{o}lder continuous on $\Omega=[-a,a]$.
In the present context, the significance of the H\"{o}lder spaces are given by the following two theorems.

\begin{theorem}\label{theo:extensionofh}  \rm
Let $h(z)$ be a Herglotz function as defined in \eref{eq:Herglotz2} with measure $\beta$,
and let $\Omega\subset\R$ be a finite union of closed and bounded intervals on the real axis and ${\cal O}$ some  neighborhood of  $\Omega$. 
Suppose that the measure $\beta$ is absolutely continuous on $\overline{\cal O}$ with
density $\beta_2^\prime\in C^{0,\alpha}(\overline{\cal O})$, $0<\alpha<1$, where $\overline{\cal O}$ 
denotes the closure of ${\cal O}$.
The measure $\beta$ is here divided into two parts, $\beta_1$ and $\beta_2$, such that $\beta_1$ is supported on $\R\setminus{\cal O}$
and $\beta_2$ is supported on  $\overline{\cal O}$.
Then the Herglotz function can be H\"{o}lder continuously extended (with H\"{o}lder exponent $\alpha$) from $\C^+$ to $\C^+\cup\Omega$ with boundary values
\begin{equation}\label{eq:HTdef}
h(x)=b_1x+c+\dashint_{\R}\left(\frac{1}{\xi-x}-\frac{\xi}{1+{\xi}^2}\right){\rm d}\beta(\xi)+\iu \pi \beta_2^\prime(x), \quad x\in\Omega,
\end{equation}
and where the integral is taken as a Cauchy principal value. 
\hfill $\blacksquare$
\end{theorem}
{\bf Proof:} {This theorem can readily be proven by employing the Sokhotski-Plemelj theorem \cite[Theorem 7.6, p.~101]{Kress1999}, or the Plemelj-Privalov theorem \cite[Theorem 5.7.21, p.~484]{Simon2015}
which are formulated for a bounded and simply connected domain in the complex plane.}
Without loss of generality, it is assumed here that the approximation domain $\Omega$ consists of a single closed and bounded interval on $\R$.
Consider the integral representation \eref{eq:Herglotz2} subjected to the limit $\Im\{z\}\rightarrow 0+$, where $\Re\{z\}\in\Omega\subset\overline{\cal O}$.

For the terms $b_1z+c$ as well as the non-singular part of the integral over $\R\setminus\overline{\cal O}$ the statement {in \eref{eq:HTdef}} is clear. For the remaining singular integral  the interval $\overline{\cal O}$ {is extended} to a closed simple smooth curve $\overline{\cal O}\cup\Gamma$, which is the boundary of a simply connected domain $D$ in $\C^+$, that is $\partial D=\overline{\cal O}\cup \Gamma$. Furthermore,   the 
density $\beta_2^\prime$  {is extended} as a H\"{o}lder continuous function to the whole of $\partial D$. 
Consider now the  function 
\begin{equation}\label{eq:SokhotskiPlemelj0}
G(z)=\int_{\partial D}\left(\frac{1}{u-z}-\frac{u}{1+u^2}\right)\beta_2^\prime(u){\rm d}u,  
\end{equation}
which is analytic in $D$. By \cite[Theorem 7.6, p.~101]{Kress1999} the function { $G(z)$} can be  
H\"{o}lder continuously {extended} (with H\"{o}lder exponent $\alpha$) to the boundary $\partial D$.
Note, that if the integral over $\partial D$ is decomposed into {integrals over} $\overline{\cal O}$ and $\Gamma$ then for $z\in \overline{\cal O}$ the integral over $\Gamma$ is analytic in $z$. Hence, for the extension of $G$ the explicit formula in \cite[Theorem 7.6, p.~101]{Kress1999} yields
\begin{eqnarray}\label{eq:HTdef2}
\lim_{\Im\{z\}\rightarrow 0+}\int_{\overline{\cal O}}\left(\frac{1}{\xi-z}-\frac{\xi}{1+{\xi}^2}\right)\beta_2^\prime(\xi){\rm d}\xi \nonumber \\
\quad =\dashint_{\overline{\cal O}}\left(\frac{1}{\xi-x}-\frac{\xi}{1+{\xi}^2}\right)\beta_2^\prime(\xi){\rm d}\xi+\iu \pi \beta_2^\prime(x), \quad x\in\Omega, 
\end{eqnarray}
and where the integral is taken as a Cauchy principal value. 
{The statement in \eref{eq:HTdef} is thus finally established.}
\hfill $\square$

\begin{theorem}\label{theo:holdercontinuityHilbert}  \rm
Let $\Omega_2\subset\Omega_1$ be finite unions of closed and bounded intervals on $\R$.
Define the Cauchy integral operator $H$ by
\begin{equation}\label{eq:HTholdercontinuitydef}
Hf(x)=\dashint_{\Omega_1}\left(\frac{1}{\xi-x}-\frac{\xi}{1+{\xi}^2}\right)f(\xi){\rm d}\xi, \quad x\in\Omega_2,
\end{equation}
where the integral is taken as a Cauchy principal value. 
Then $H:C^{0,\alpha}(\Omega_1)\rightarrow C^{0,\alpha}(\Omega_2)$ is a bounded operator.
\hfill $\blacksquare$
\end{theorem}

{\bf Proof:} To prove this theorem, it is noted that the singular part of the operator in \eref{eq:HTholdercontinuitydef}
can be treated in the same way as the Cauchy integral operator in \cite[Theorem 7.6 and Corollary 7.7 on pp.~101-102]{Kress1999}.
See also the proof of Theorem \ref{theo:extensionofh} above.
\hfill $\square$

\section{Passive approximation}\label{sect:passiveapprox}

\subsection{On the density of finite B-spline expansions in H\"{o}lder space}
Let $\Omega\subset\R$ denote the approximation domain defined as a finite union of closed and bounded intervals on the real axis.
Finite-dimensional approximations of continuous functions on $\Omega$ can be generated by {\em finite spline functions}.
Here, a finite spline function on $\Omega$ of order $m\geq 2$ is an $m-2$ times continuously differentiable function which is piecewise polynomial of order at most $m-1$.
The finitely many break-points of a  spline function are referred to as knots or partition of  $\Omega$, as described in \eg \cite{Boor1968,Boor1972,Dahlquist+Bjork1974}. 

When the distance between any two consecutive knots are equal the partition is said to be {\em uniform} 
and the corresponding spline function is a {\em uniform spline function}.
Note that any spline function is H\"{o}lder continuous with H\"{o}lder exponent in $(0,1)$.

A B-spline of order $m\geq 2$ is a compactly supported positive basis spline function which is 
piecewise polynomial of order $m-1$, \ie linear, quadratic, cubic, etc.,
and which is defined by $m+1$ break-points as described in \eg \cite{Boor1972}. A single B-spline extends over $m+1$ knots and hence the {\em internal} knots belonging to $\Omega$ must be extended with at least  
$m-1$ {\em external} knots on each side of an interval of $\Omega$ to give full support to the first and the last B-spline on $\Omega$. 

The set of all finite spline functions on $\Omega$ of order $m$ 
is a linear space.
Given a  function $f$ of this space
there exist finitely many B-splines of order $m$ such that $f$ can be written as a linear combination of these B-splines and, moreover, this representation is unique, \cf \cite{Boor1972}. Hence, the set of all B-splines on $\Omega$ of order $m$ constitutes a (well-conditioned) basis of this linear space, which thus
 will also be referred to here as the linear space of finite $m$th order B-spline expansions on $\Omega$. 
Note that the partitions can be chosen arbitrarily and the space of finite $m$th order B-spline expansions on $\Omega$ is therefore an infinite-dimensional subspace of $C^{0,\alpha}(\Omega)$ where $0<\alpha<1$.

The proof of the Theorem \ref{theo:BsplinesinHolderspace} given below is constructive and therefore restricted to the linear space of finite uniform $m$th order B-spline expansions on $\Omega$.
Without loss of generality, a single interval $\Omega=[0,|\Omega|]$ is considered here with the following partition and internal knots 
\begin{equation}\label{eq:partitiondef}
\left\{\begin{array}{ll}
\displaystyle x_k=k\frac{|\Omega|}{N} & k=0,\ldots,N, \vspace{0.2cm} \\
\displaystyle I_k=[x_{k-1},x_k) & k=1,\ldots,N,  \vspace{0.2cm} \\
\end{array}\right.
\end{equation}
where $N$ is the number of subintervals $I_k$ in $\Omega$. This partition is also associated with $m-1$ uniformly spaced external knots $x_k$ 
($k=1-m,\ldots,-1$ and $k=N+1,\ldots,N+m-1$) with corresponding subintervals $I_k$ on each side of the interval $\Omega$.

\begin{figure}[htb!]
\begin{center}
\includegraphics[width=0.95\textwidth]{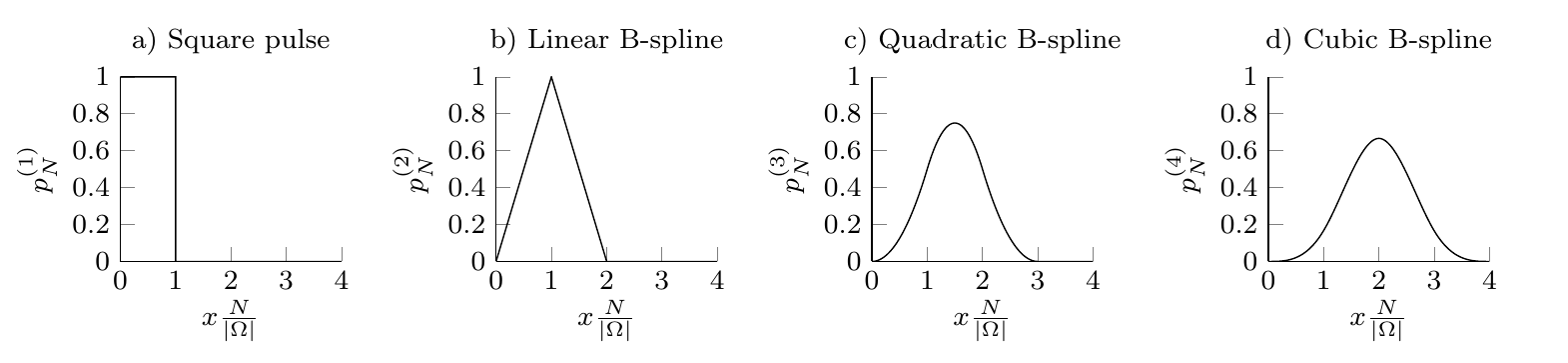}
\end{center}
\vspace{-0.5cm}
\caption{Basic square pulse $p_N^{(1)}$ and the generated prototype B-spline basis functions $p_N^{(m)}=g_N^{(1)}*p_N^{(m-1)}$ of order $m=2,3,4$, \ie linear, quadratic and cubic, respectively.}
\label{fig:matfig1}
\end{figure}

Consider first the $m$th order B-spline $p_N^{(m)}$ which is supported on 
$K_N^{(m)}=[0,\frac{m|\Omega|}{N}]$, the other basis functions can then be generated as translations of this prototype. It is 
defined here by the second order (linear) prototype B-spline $p_N^{(2)}$ as
\begin{equation}\label{eq:pNmdef}
p_N^{(m)}=g_N^{(m-2)}*p_N^{(2)}, \quad \mrm{where} \quad p_N^{(2)}(x)
=\left\{\begin{array}{ll}
\frac{N}{|\Omega|}x & 0\leq x \leq \frac{|\Omega|}{N},  \\
2-\frac{N}{|\Omega|}x & \frac{|\Omega|}{N} \leq x \leq \frac{2|\Omega|}{N},  \\
0 & \mrm{elsewhere},
\end{array}\right.
\end{equation}
and where $g_N^{(0)}=\delta$ gives an identity operation and 
$g_N^{(l)}$ for $l\geq 1$ is an $l$ times repeated convolution with the normalized square pulse
\begin{equation}
g_N^{(1)}(x)=\left\{\begin{array}{ll}
\frac{N}{|\Omega|} & 0\leq x \leq \frac{|\Omega|}{N},  \\
0 & \mrm{elsewhere}.
\end{array}\right.
\end{equation}
This definition is equivalent to the recursion
$p_N^{(m)}=g_N^{(1)}*p_N^{(m-1)}$ for $m=2,3,\ldots$, where $p_N^{(1)}$ is a
square pulse as illustrated in Figure~\ref{fig:matfig1}.
Note that the generic convolution operand $g_N^{(l)}$ in \eref{eq:pNmdef} is normalized to that $\int_{\R}g_N^{(l)}(x)\mrm{d}x=1$.
For any given $f\in C^{0,\alpha}(\Omega)$, a uniform 2nd order (linear) B-spline approximation is given by
\begin{equation}\label{eq:fNdef}
f_N(x)=\sum_{k=2-m}^{N}f_k p_N^{(2)}(x-x_{k-1}), \quad x\in\R,
\end{equation}
where $f_k=f(x_k)$ and where $f$ has been H\"{o}lder continuously extended for $x<0$ (\eg $f(x)=f(0)$ for $x<0$) to 
provide for the generation of higher order B-splines {having full support also outside the interval $\Omega$}.
A general uniform $m$th order B-spline approximation of $f$ can now be obtained as 
\begin{equation}\label{eq:gNconvfN}
g_N^{(m-2)}(x)*f_N(x)=\sum_{k=2-m}^{N}f_k p_N^{(m)}(x-x_{k-1}), \quad x\in\R.
\end{equation}
Note that the function $g_N^{(m-2)}$ (interpreted as a distribution) will approach the Dirac distribution as $N\rightarrow \infty$.
Hence, if $f_N$ converges to $f$ in some suitable sense it is also expected that $g_N^{(m-2)}*f_N$ will converge to $f$ as $N\rightarrow \infty$.
A precise notion for the density of these B-spline expansions in $C^{0,\alpha}(\Omega)$ is given below.

\begin{theorem}\label{theo:BsplinesinHolderspace}  \rm
Let $0<\alpha<1$, $m\geq 2$ and let $\Omega$ be a finite union of closed and bounded intervals on $\R$.
The linear space of finite uniform $m$th order B-spline expansions on $\Omega$ is dense in $C^{0,\alpha}(\Omega)$
in the topology induced by any of the $C^{0,\alpha^\prime}(\Omega)$-norms
where $0<\alpha^\prime<\alpha$.  \hfill $\blacksquare$ 
\end{theorem}

{\bf Proof:} Without loss of generality, the interval $\Omega=[0,|\Omega|]$ is considered with a partition consisting of $N$ subintervals and approximating functions
defined as in \eref{eq:partitiondef} through \eref{eq:gNconvfN} above. 
Let $f\in C^{0,\alpha}(\Omega)$ where $0<\alpha<1$, so that $|f(x)-f(y)|\leq C|x-y|^\alpha$ for all $x,y\in\Omega$ and let $0<\alpha^\prime<\alpha$.
Due to the uniform continuity of $f$, it is readily seen that $g_N^{(m-2)}*f_N\rightarrow f$ in max-norm.
Hence, only the second part of the definition \eref{eq:Holdernorm} is treated in detail, and the task is to show that
\begin{equation}
\sup_{x,y\in\Omega,x\neq y} \frac{\left| f(x)-g_N^{(m-2)}(x)*f_N(x) - (f(y)-g_N^{(m-2)}(y)*f_N(y)) \right| }{|x-y|^{\alpha^\prime}}
\end{equation}
converges to zero as $N\rightarrow\infty$.

Consider first the case when $|x-y|<\frac{|\Omega|}{N}$ and let $x<y$. The following bounds are obtained
\begin{eqnarray}
\left| f(x)-g_N^{(m-2)}(x)*f_N(x) - (f(y)-g_N^{(m-2)}(y)*f_N(y)) \right| \nonumber \\
\quad \leq\left|f(x)-f(y)\right| + \left|g_N^{(m-2)}(x)*f_N(x)-g_N^{(m-2)}(y)*f_N(y)\right| \nonumber \\
\quad=\left|f(x)-f(y)\right| + \left|\int_{K_N^{(m-2)}}g_N^{(m-2)}(\xi)\left(f_N(x-\xi)-f_N(y-\xi) \right)\mrm{d}\xi\right| \nonumber \\
\quad\leq\left|f(x)-f(y)\right| + \max_{\xi\in K_N^{(m-2)}}\left|f_N(x-\xi)-f_N(y-\xi)\right|.
\end{eqnarray}
Write
\begin{equation}\label{eq:fNdefs}
\left\{\begin{array}{ll}
\displaystyle f_N(x)=\frac{N}{|\Omega|}\left(f_k-f_{k-1} \right)(x-x_k)+f_k, & x\in I_k, \vspace{0.2cm} \\
\displaystyle f_N(y)=\frac{N}{|\Omega|}\left(f_{k+1}-f_{k} \right)(y-x_k)+f_k, & y\in I_{k+1},
\end{array}\right.
\end{equation}
and consider the following bounds 
\begin{eqnarray}
\left|f_N(x-\xi)-f_N(y-\xi)\right| \nonumber \\
\quad\leq\frac{N}{|\Omega|}\left|f_k-f_{k-1}\right|\left| x-y\right|  \leq \frac{N}{|\Omega|}C\left(\frac{|\Omega|}{N}\right)^\alpha |x-y|,
\end{eqnarray}
where $x-\xi\in I_k$, $y-\xi\in I_k$ and $\xi\in K_N^{(m-2)}$, and
\begin{eqnarray}
\left|f_N(x-\xi)-f_N(y-\xi)\right| \nonumber \\
\quad \leq\frac{N}{|\Omega|}\left( |f_k-f_{k-1}| |x-\xi-x_k|+|f_{k+1}-f_{k}| |y-\xi-x_k|\right)   \nonumber \\
\quad \leq 2\frac{N}{|\Omega|}C\left(\frac{|\Omega|}{N}\right)^\alpha|x-y|,
\end{eqnarray}
where $x-\xi\in I_k$, $y-\xi\in I_{k+1}$ and $\xi\in K_N^{(m-2)}$. In conclusion, 
\begin{equation}
\max_{\xi\in K_N^{(m-2)}}\left|f_N(x-\xi)-f_N(y-\xi)\right|\leq 2C\left(\frac{|\Omega|}{N}\right)^{\alpha-1} |x-y|,
\end{equation}
and finally
\begin{eqnarray}\label{eq:mBsplineconvterm1}
\sup_{|x-y|<\frac{|\Omega|}{N},x\neq y} \frac{\left| f(x)-g_N^{(m-2)}(x)*f_N(x) - (f(y)-g_N^{(m-2)}(y)*f_N(y)) \right| }{|x-y|^{\alpha^\prime}} \nonumber \\
\quad \leq \sup_{|x-y|<\frac{|\Omega|}{N},x\neq y} \frac{\left|f(x)-f(y)\right| + \max_{\xi\in K_N^{(m-2)}}\left|f_N(x-\xi)-f_N(y-\xi)\right| }{|x-y|^{\alpha^\prime}} \nonumber \\
\quad \leq \sup_{|x-y|<\frac{|\Omega|}{N},x\neq y}\frac{C|x-y|^\alpha+2C\left(\frac{|\Omega|}{N}\right)^{\alpha-1} |x-y|}{|x-y|^{\alpha^\prime}} \nonumber \\
\quad \leq \sup_{|x-y|<\frac{|\Omega|}{N},x\neq y}C|x-y|^{\alpha-\alpha^\prime}+2C\left(\frac{|\Omega|}{N}\right)^{\alpha-1} |x-y|^{1-\alpha^\prime} \nonumber \\
 \quad =C\left(\frac{|\Omega|}{N}\right)^{\alpha-\alpha^\prime}+2C\left(\frac{|\Omega|}{N}\right)^{\alpha-1} \left(\frac{|\Omega|}{N}\right)^{1-\alpha^\prime}
 =3C\left(\frac{|\Omega|}{N}\right)^{\alpha-\alpha^\prime}.
\end{eqnarray}

Next, consider the case when $|x-y|\geq\frac{|\Omega|}{N}$. The following bounds are obtained
\begin{eqnarray}\label{eq:upperboundsgeqON}
\left| f(x)-g_N^{(m-2)}(x)*f_N(x) - (f(y)-g_N^{(m-2)}(y)*f_N(y)) \right|   \nonumber \\
\leq \left| \int_{K_N^{(m-2)}}g_N^{(m-2)}(\xi)(f(x)-f_N(x-\xi))\mrm{d}\xi \right| \nonumber \\
\quad \quad \quad +\left| \int_{K_N^{(m-2)}}g_N^{(m-2)}(\eta)(f(y)-f_N(y-\eta))\mrm{d}\eta \right|\nonumber \\
 \leq  \max_{\xi\in K_N^{(m-2)}}|f(x)-f_N(x-\xi)| +\max_{\eta\in K_N^{(m-2)}}|f(y)-f_N(y-\eta)|.
\end{eqnarray}
Employ the first expression in \eref{eq:fNdefs} to obtain the following bounds
\begin{eqnarray}
|f(x)-f_N(x-\xi)|\leq |f(x)-f_{k}|+\frac{N}{|\Omega|}|f_{k}-f_{k-1}| |x-\xi-x_{k} |\nonumber \\
 \leq C\left( \frac{|\Omega|}{N}(m-1)\right)^\alpha+\frac{N}{|\Omega|}C\left( \frac{|\Omega|}{N}\right)^\alpha \frac{|\Omega|}{N}\nonumber \\
\quad =C\left( \frac{|\Omega|}{N}\right)^\alpha\left(1+(m-1)^\alpha \right),
\end{eqnarray}
where $x-\xi\in I_k$ and $\xi\in K_N^{(m-2)}$, and similarly for the second term in \eref{eq:upperboundsgeqON}.
Finally,
\begin{eqnarray}\label{eq:mBsplineconvterm2}
\sup_{|x-y|\geq\frac{|\Omega|}{N}}\frac{\left| f(x)-g_N^{(m-2)}(x)*f_N(x) - (f(y)-g_N^{(m-2)}(y)*f_N(y)) \right| }{|x-y|^{\alpha^\prime}} \nonumber \\
\leq \sup_{|x-y|\geq\frac{|\Omega|}{N}}\frac{\max_{\xi\in K_N^{(m-2)}}|f(x)-f_N(x-\xi)|+\max_{\eta\in K_N^{(m-2)}}|f(y)-f_N(y-\eta)|}{|x-y|^{\alpha^\prime}}  \nonumber \\
\leq \sup_{|x-y|\geq\frac{|\Omega|}{N}} \frac{2C\left( \frac{|\Omega|}{N}\right)^\alpha\left(1+(m-1)^\alpha \right)}{|x-y|^{\alpha^\prime}}  \nonumber \\
\quad \quad = 2C\left( \frac{|\Omega|}{N}\right)^{\alpha-\alpha^\prime}\left(1+(m-1)^\alpha\right).
\end{eqnarray}
The right hand side in the two cases \eref{eq:mBsplineconvterm1} and \eref{eq:mBsplineconvterm2}
all converge to zero as $N\rightarrow\infty$, which establishes the theorem.
Note also that the case $m=2$ can be treated similarly but without involving the convolution integral.
\hfill $\square$

\subsection{Topology for Herglotz functions that are locally H\"{o}lder continuous on the real axis}
Let $\Omega$ be defined as above and let $C(\Omega)$ denote the Banach space consisting of all complex valued continuous functions defined on $\Omega$, which is equipped with the usual max-norm denoted $\|\cdot\|_\infty$.
Further, let $w>0$ denote a positive continuous weight function on $\Omega$ and let $\mrm{L}^p(w,\Omega)$ denote the Banach space with norm 
\begin{equation}
\|f\|_{\mrm{L}^p(w,\Omega)}=\left(\int_\Omega w(x)|f(x)|^p{\rm d}x\right)^{1/p},
\end{equation}
where $1\leq p<\infty$, \cf \cite{Rudin1987}. Similarly, let $\mrm{L}^\infty(w,\Omega)$ denote the Banach space with the norm $\|f\|_{\mrm{L}^\infty(w,\Omega)}$
defined by taking the essential supremum \cite{Rudin1987} of the function $w|f|$.
The following inclusions are now valid $C^{0,\alpha}(\Omega)\subset C^{0,\alpha^\prime}(\Omega)\subset C(\Omega)\subset \mrm{L}^p(w,\Omega)$,
where $0<\alpha^\prime<\alpha<1$ and $1\leq p\leq\infty$.

\begin{definition} \label{def:coneV} \rm
Let $V^{\alpha,p}(w,\Omega)\subset \mrm{L}^p(w,\Omega)$, $0<\alpha<1$, $1\leq p\leq \infty$, 
denote the convex cone consisting of all complex valued functions $h\in \mrm{L}^p(w,\Omega)$ 
such that $h$ is the H\"{o}lder continuous extension (with exponent $\alpha$) to $\overline{\cal O}$ of a Herglotz function where $\overline{\cal O}$ is the closure of some neighborhood ${\cal O}$ of $\Omega$.
Note that the set $V^{\alpha,p}(w,\Omega)$ is independent of $p$ and $w$ but is equipped with the topology of $\mrm{L}^p(w,\Omega)$.
\hfill $\square$
\end{definition}

\begin{definition} \label{def:coneVspline} \rm
Let $V_{m}^{p}(w,\Omega)\subset V^{\alpha,p}(w,\Omega)$, $m\geq 2$, $1\leq p\leq \infty$, $0<\alpha<1$,
denote the convex cone consisting of all complex valued continuous functions $h\in V^{\alpha,p}(w,\Omega)$ 
such that $h$ is the H\"{o}lder continuous extension to $\Omega$ of a Herglotz function, where the corresponding measure $\beta$ is absolutely continuous on $\R$ and
where the density $\beta^\prime$ is a finite $m$th order B-spline expansion.
\hfill $\square$
\end{definition}

The following topological result can now be established.
\begin{theorem}\label{theo:densesetspline}  \rm
 The convex cone $V_{m}^p(w,\Omega)$ with $m\geq 2$ is dense in $V^{\alpha,p}(w,\Omega)$ in the 
 topology induced by the norm on $\mrm{L}^p(w,\Omega)$, and where $0<\alpha<1$, $1\leq p \leq \infty$. 
  \hfill $\blacksquare$
\end{theorem}

{\bf Proof:} Without loss of generality, and for simplicity, a single interval $\Omega$ is considered below.
Let  $h\in V^{\alpha,p}(w,\Omega)$ be given. 
It will be shown that for every  $\varepsilon>0$ there exists a function 
$\widetilde h\in V_{m}^p(w,\Omega)$ such that $\|\widetilde h-h\|_{\mrm{L}^p(w,\Omega)}<\varepsilon$.
Since $h\in V^{\alpha,p}(w,\Omega)$, there exists an open set $\cal O\subset\R$ such that 
according to Definition \ref{def:coneV} and Theorems \ref{theo:limitsofbetaandimh} and \ref{theo:extensionofh} the function $h$ can be represented by
\begin{eqnarray}\label{eq:hdef}
h(x)=b_1 x+c+ \int_{\R\setminus\overline{\cal O}}\left(\frac{1}{\xi-x}-\frac{\xi}{1+{\xi}^2}\right){\rm d}\beta_1(\xi) \\[1mm]
\hspace*{28mm} +\dashint_{\overline{\cal O}}\left(\frac{1}{\xi-x}-\frac{\xi}{1+{\xi}^2}\right)\beta_2^\prime(\xi){\rm d}\xi +\iu\pi\beta_2^\prime(x), \nonumber
\end{eqnarray}
for $x\in\Omega$. Here, the measure $\beta$ is divided into the two parts, $\beta_1$ and $\beta_2$, supported on $\R\setminus{\cal O}$ and $\overline{\cal O}$, respectively, and where $\beta_2$ is absolutely continuous on $\overline{\cal O}$ with density $\beta_2^\prime\in C^{0,\alpha}(\overline{\cal O})$.
Without loss of generality, it is assumed that  $b_1=0$ and $c=0$ in the representation \eref{eq:hdef} of $h$.

Furthermore, given the set $\cal O$,  any function $\widetilde h\in V_{m}^p(w,\Omega)$ with $\widetilde b_1=\widetilde c=0$ can
be  represented similarly as
\begin{eqnarray}\label{eq:hndef}
\widetilde h(x)=\int_{\R\setminus\overline{\cal O}}\left(\frac{1}{\xi-x}-\frac{\xi}{1+{\xi}^2}\right){\rm d}\widetilde\beta_{1}(\xi) \\[1mm]
\hspace*{18mm} +\dashint_{\overline{\cal O}}\left(\frac{1}{\xi-x}-\frac{\xi}{1+{\xi}^2}\right) \widetilde\beta_{2}^\prime(\xi){\rm d}\xi 
+\iu\pi\widetilde \beta_{2}^\prime(x), \nonumber
\end{eqnarray}
for $x\in\Omega\subset\overline{\cal O}\subset\R$, and where the corresponding measure $\widetilde \beta$ is divided into the two parts,
$\widetilde \beta_{1}$ and $\widetilde \beta_{2}$, supported on $\R\setminus{\cal O}$ and $\overline{\cal O}$, respectively.
Note that both $\widetilde \beta_{1}$ and $\widetilde \beta_{2}$ have compact support, and are absolutely continuous with   
densities $\widetilde \beta_{1}^\prime(x)$ and $\widetilde \beta_{2}^\prime(x)$, respectively, which are  positive, continuous and piecewise polynomial functions.

Let $r>0$  be sufficiently large such that for  $R_r=(-r,r)\subset\R$ it holds that $R_r\supset \overline{\cal O}$. 
By furthermore restricting the measure $\widetilde \beta_1$ such that  $R_r\supset\textrm{supp}\{\widetilde \beta_{1}\}$ 
the approximation error $\|\widetilde h-h\|_{\mrm{L}^p(w,\Omega)}$ can be estimated as
\begin{eqnarray}
\hspace*{-2cm}\|\widetilde h-h\|_{\mrm{L}^p(w,\Omega)} & \leq &  \left\|\int_{\R\setminus R_{r}}\left(\frac{1+\xi x}{\xi-x}\right)\frac{{\rm d}\beta_1(\xi)}{1+\xi^2}\right\|_{\mrm{L}^p(w,\Omega)}  \label{eq:term1} \\[1mm]
&& + \left\|\int_{R_{r}\setminus \overline{\cal O}}\left(\frac{1}{\xi-x}-\frac{\xi}{1+{\xi}^2}\right)\left({\rm d}\widetilde\beta_{1}(\xi)-{\rm d}\beta_1(\xi)\right) \right\|_{\mrm{L}^p(w,\Omega)}  \label{eq:term2} \\[1mm]
&& + \left\|\dashint_{\overline{\cal O}}\left(\frac{1}{\xi-x}-\frac{\xi}{1+{\xi}^2}\right)\left(\widetilde\beta_{2}^\prime(\xi)-\beta_2^\prime(\xi)\right){\rm d}\xi \right\|_{\mrm{L}^p(w,\Omega)}  \label{eq:term3} \\[1mm]
&& + \left\|\iu\pi\left(\widetilde\beta_{2}^\prime(x)-\beta_2^\prime(x)\right)\right\|_{\mrm{L}^p(w,\Omega)}. \label{eq:term4}
\end{eqnarray}

\begin{enumerate}
\item For the term in \eref{eq:term1},
the convergence property of the measure $\int_{\R}{\rm d}\beta_1(\xi)/(1+\xi^2)<\infty$  implies that
this term can be made  arbitrarily small by choosing $r$ sufficiently large.

\item Next, it will be shown that there exists a finite spline function $\widetilde \beta^\prime_1$ such that 
the term in \eref{eq:term2} converges to zero. To this end, $\beta_1$ will first be approximated with an arbitrary continuous function and then by a finite spline function. 

To start with,  choose a cutoff function  $\psi(\xi)$, \ie a $C^\infty(\R)$-function with compact support, and  such that $\psi(\xi)\equiv 1$ in a neighborhood of
$I=R_{r}\setminus \overline{\cal O}$, \cf \cite{Hormander1983}.
Hence, also $\phi_x(\xi)=\big(\frac{1}{\xi-x}-\frac{\xi}{1+{\xi}^2}\big)\psi(\xi)$ is a test function with compact support.
Then the integral in \eref{eq:term2} can be rewritten using the notion of distributions as
$$
\int_{R_{r}\setminus \overline{\cal O}}\left(\frac{1}{\xi-x}-\frac{\xi}{1+{\xi}^2}\right)\left(\widetilde \beta^\prime_1(\xi){\rm d}\xi-{\rm d}\beta_1(\xi)\right) = \langle \widetilde\beta_{1}^{\rm }-\beta_{1},\phi_x(\xi) \rangle,
$$ 
where here and in the following the same symbol is used for a measure and the corresponding zero-order distribution.
Now, choose a  continuous function $g(\xi)\in C^0(\R)$  with  $g(\xi)\geq 0$ and  $\int_{K}g(\xi){\rm d}\xi=1$ which has compact support and denote  $K=\textrm{supp}\{g\}\subset\R$. 
Set $g_n(\xi)=ng(n\xi)$ and note that $\lim_{n\rightarrow\infty}\textrm{supp}\{g_n(\xi)\}=\{0\}$, 
\cf \cite[Theorems 4.1.4-5 on p.~89]{Hormander1983}.

Now,  define the distribution $\beta_{1n}$ by the convolution $\beta_{1n}=\beta_1*g_n(\xi)$, and note that it  corresponds to a continuous function $\beta_{1n}^{\prime}(\xi)\in C^0(\R)$.
The following distributional (convolution) identity is now obtained
\begin{equation}\label{eq:distrrelgnphix1}
\langle \beta_{1n}-\beta_{1},\phi_x(\xi) \rangle=\langle  \beta_{1},g_n(\xi)*\phi_x(\xi)-\phi_x(\xi) \rangle.
\end{equation}
The right hand side can be explicitly calculated as
\begin{eqnarray}\label{eq:distrrelgnphix2}
\int_{I}\left[\int_{\R}ng(nu)\phi_x(\xi-u){\rm d}u-\phi_x(\xi) \right]{{\rm d}\beta_1(\xi)} \nonumber \\
\quad =\int_{I}\left[\int_{K}g(\tau)\left[\phi_x(\xi-\tau/n)-\phi_x(\xi)\right]{\rm d}\tau \right]{{\rm d}\beta_1(\xi)}\nonumber  \\
 \quad  =\int_{I}\left[\int_{K}g(\tau)\left(\frac{\tau/n}{(\xi-x)(\xi-\tau/n-x)}+\frac{r_n(\xi,\tau)}n\right){\rm d}\tau \right]{{\rm d}\beta_1(\xi)},
\end{eqnarray}
where $r_n$ is a continuous function on $I\times K$ and uniformly bounded in $n$.  
Here, it has been exploited that $\psi(\xi)\equiv 1$ on $I$ and $\psi(\xi-\tau/n)\equiv 1$ on $I$ when $n$ is sufficiently large.
In order to estimate the  integral in \eref{eq:distrrelgnphix2}, it is noted that   $I=R_{r}\setminus \overline{\cal O}$   consists of two intervals. 
Without loss of generality the right one is considered here, i.e.  $ x <\xi$, and hence there exist points $x_0$ and $x_1$ such that   $x_0\leq x \leq x_1<\xi$. 
Then  \eref{eq:distrrelgnphix1} and 
\eref{eq:distrrelgnphix2} yield for all $x\in\Omega$
\begin{eqnarray}\label{eq:distrrelgnphix3}
\hspace*{-22mm}\left| \langle \beta_{1n}-\beta_{1},\phi_x(\xi) \rangle \right|
=\left|\int_{I}\left[\int_{K}g(\tau)\left(\frac{\tau/n}{(\xi-x)(\xi-\tau/n-x)}+\frac{r_n(\xi,\tau)}n\right){\rm d}\tau \right]{{\rm d}\beta_1(\xi)}\right|\nonumber \\
 \hspace*{-3mm}\leq \int_{I}\left[\int_{K}g(\tau)\left(\frac{|\tau/n|}{(\xi-x_1)(\xi-\tau/n-x_1)}+\frac{r_n(\xi,\tau)}n\right){\rm d}\tau \right]{{\rm d}\beta_1(\xi)}.
\end{eqnarray}
For the left interval a similar bound can be found and hence, $n$ can be chosen sufficiently large to make \eref{eq:distrrelgnphix3} arbitrarily small uniformly over $x\in\Omega$.
Next,  for any given $n$, let  $\widetilde \beta_{1}$ be an arbitrary $m$th order B-spline approximation of  $\beta_{1n}$ on $R_{r}\setminus \overline{\cal O}$. 
Then 
\begin{equation}\label{eq:Bsplineestimate}
\left|\langle \widetilde\beta_{1}-\beta_{1n}, \phi_x(\xi)\rangle\right|
\leq \sup_{\xi\in R_{r}\setminus \overline{\cal O}}\left|\widetilde\beta_{1}^{\prime}(\xi)-\beta_{1n}^{\prime}(\xi)\right|M\int_{R_{r}\setminus \overline{\cal O}}{\rm d}\xi,
\end{equation}
where $M$ is the supremum of $\left|\phi_x(\xi)\right|$ over $(\xi,x)\in I\times\Omega$.
Since the finite $m$th order B-spline expansions are dense in the space of continuous functions with compact support $C_0^0(\R)$ 
in the usual max-norm {(can be shown similarly as in the proof of Theorem \ref{theo:BsplinesinHolderspace} above)}, 
it follows that $\widetilde \beta_{1}^\prime$ can be chosen such that also the left hand side of \eref{eq:Bsplineestimate} becomes 
arbitrarily small uniformly over $x\in\Omega$.
Finally, by employing the triangle inequality 
\begin{equation}\label{eq:triangleineq}
\hspace*{-15mm}\left|\langle \widetilde\beta_{1}-\beta_{1}, \phi_x(\xi)\rangle\right|
\leq \left|\langle \widetilde\beta_{1}-\beta_{1n}, \phi_x(\xi)\rangle\right|
+\left|\langle \beta_{1n}-\beta_{1}, \phi_x(\xi)\rangle\right|,
\end{equation}
it is now seen that the left hand side of \eref{eq:triangleineq} can be made arbitrarily small uniformly over $x\in\Omega$.
The convergence of the term \eref{eq:term2} is finally established by using the inclusion $C(\Omega)\subset \mrm{L}^p(w,\Omega)$
and the inequalities $\|f\|_{\mrm{L}^p(w,\Omega)}\leq \|f\|_\infty\left(\int_\Omega w(x)\mrm{d}x \right)^{1/p}$ for $1\leq p <\infty$
and $\|f\|_{\mrm{L}^\infty(w,\Omega)}\leq \|f\|_\infty\|w\|_\infty$.

\item For the final terms in \eref{eq:term3} and \eref{eq:term4}, it is not possible to employ convergence in $C(\Omega)$
since the Cauchy integral operator is an unbounded operator in this space. Instead, the Theorems \ref{theo:holdercontinuityHilbert} and \ref{theo:BsplinesinHolderspace} are used based on convergence in H\"{o}lder space.
In particular, it is assumed that $\beta_2^\prime\in C^{0,\alpha}(\overline{\cal O})$ where $0<\alpha<1$. 
Now, fix some $\alpha^\prime$ with $0<\alpha^\prime<\alpha$ and note that $C^{0,\alpha}(\overline{\cal O})\subset C^{0,\alpha^\prime}(\overline{\cal O})\subset C(\Omega)\subset \mrm{L}^p(w,\Omega)$. 
It then follows from Theorem \ref{theo:BsplinesinHolderspace} that there exists a sequence of finite uniform $m$th order B-spline expansions $\widetilde \beta^\prime_{2n}$ on $\Omega$
such that $\|\widetilde \beta^\prime_{2n}-\beta_2^\prime\|_{\alpha^\prime}\rightarrow 0$ as $n\rightarrow\infty$. It follows immediately that
$\|\widetilde \beta^\prime_{2n}-\beta_2^\prime\|_{\mrm{L}^p(w,\Omega)}\rightarrow 0$ as $n\rightarrow\infty$, which establishes the convergence of the term \eref{eq:term4}.
Consider now the Cauchy integral operator $H: C^{0,\alpha^\prime}(\overline{\cal O})\to C^{0,\alpha^\prime}(\Omega)$ as defined in \eref{eq:HTholdercontinuitydef} and which is
a bounded operator by Theorem \ref{theo:holdercontinuityHilbert}. Since  $\|\widetilde \beta^\prime_{2n}-\beta_2^\prime\|_{\alpha^\prime}\rightarrow 0$ as $n\rightarrow\infty$, it follows from
the boundedness of $H$ that $\|H(\widetilde \beta^\prime_{2n}-\beta_2^\prime)\|_{\alpha^\prime}\rightarrow 0$ as $n\rightarrow\infty$ 
and finally that $\|H(\widetilde \beta^\prime_{2n}-\beta_2^\prime)\|_{\mrm{L}^p(w,\Omega)}\rightarrow 0$
as $n\rightarrow\infty$, which establishes the convergence of the term \eref{eq:term3}.

\end{enumerate}
\hfill $\square$

\subsection{The passive approximation problem}
The following passive approximation problem is considered.
\begin{definition} \label{def:approxformulation} \rm
Let $F\in C(\Omega)$ and consider the problem to approximate $F$ based on 
the set of Herglotz functions $h\in V^{\alpha,p}(w,\Omega)$ having H\"{o}lder continuous extensions to $\Omega$ for $0<\alpha<1$, $1\leq p\leq \infty$, see Definition \ref{def:coneV}.
The greatest lower bound on the approximation error over the convex cone $V^{\alpha,p}(w,\Omega)$
is defined by
 \begin{equation}\label{eq:dinfdef}
d= \displaystyle  \inf_{h\in V^{\alpha,p}(w,\Omega)}\| h-F \|_{\mrm{L}^p(w,\Omega)}.
\end{equation}
\hfill $\square$
\end{definition}
Note that the distance $d$ is independent of the H\"{o}lder exponent $\alpha$.
This follows from an application of the triangle inequality together with
the fact that $C^{0,\alpha}(\Omega)$ is dense in
$C^{0,\alpha^\prime}(\Omega)$ in the topology induced by the norm on $C(\Omega)$ (or $\mrm{L}^p(w,\Omega)$) where $0<\alpha^\prime<\alpha<1$, \cf \cite[Theorem 7.4 on pp.~96-97]{Kress1999}.

The main result of this paper can now be summarized in the following theorem.
\begin{theorem} \label{theo:cvxoptsummary} \rm
The greatest lower bound on the approximation error defined in \eref{eq:dinfdef}  is given by
 \begin{equation}\label{eq:cvxdefsummary}
d=\inf_{h\in V_m^{p}(w,\Omega)}\| h-F \|_{\mrm{L}^p(w,\Omega)},
\end{equation}
where $m\geq 2$.
\hfill $\blacksquare$
\end{theorem}
{\bf Proof:} The result in Theorem \ref{theo:cvxoptsummary} is a simple consequence of the triangle inequality together with 
the fact that $V_{m}^p(w,\Omega)$ is dense in $V^{\alpha,p}(w,\Omega)$ in the topology induced by the norm on $\mrm{L}^p(w,\Omega)$, 
and where $m\geq 2$, $0<\alpha<1$, $1\leq p \leq \infty$, see Theorem  \ref{theo:densesetspline}. 
\hfill $\square$


\begin{example}\label{ex:approxmaxnorm}\rm
In the case with the max-norm $\|\cdot\|_\infty$ (or the $\mrm{L}^\infty(w,\Omega)$-norm), there is an interesting class of problems for which there exist
non-trivial lower bounds on the distance $d$ defined in \eref{eq:dinfdef}.
The following is a simple extension of the results given in \cite{Gustafsson+Sjoberg2010a} regarding the approximation of metamaterials. 
Suppose that $F=-h_0$ is the negative of a Herglotz function $h_0$ which can be extended continuously
to $\C^+\cup\overline{\cal O}$, ${\cal O}\supset\Omega$, and which has the large argument asymptotics $h_0(z)=b_1^0z+o(z)$ as $z\toh \infty$. 
It can then be shown that the following non-trivial lower bound holds
\begin{equation}\label{eq:nontrivialbound}
\displaystyle  \| h-F \|_{\infty}\geq (b_1+b_1^0)\frac{1}{2}|\Omega|, 
\end{equation}
for all $h\in V^{\alpha,\infty}(w,\Omega)$ such that $h(z)=b_1z+o(z)$ as $z\toh \infty$, and where
$|\Omega|$ is the length of the interval $\Omega$. 
A simple example is with $F(x)=-Cx$ where $C>0$ and where $d\geq C\frac{1}{2}|\Omega|$. 
The proof of \eref{eq:nontrivialbound} is based on the theory of Herglotz functions and its associated sum rules and 
is given in \ref{sect:nontrivialbound}. 
\hfill $\square$
\end{example}

\subsection{Convex optimization with B-splines}\label{sect:cvxsect}
The impulse response of a passive system can usually be considered to be real-valued. 
In this case, the approximating Herglotz function is symmetric and can be represented as in \eref{eq:Herglotz3}.
Here, the approximating Herglotz function is furthermore assumed to be generated by a finite uniform $m$th order B-spline expansion
with extension $h\in V_{m}^p(w,\Omega)$ (see Definition \ref{def:coneVspline})
given by
\begin{equation}\label{eq:HTdefsymm}
h(x)=b_1x+\dashint_{\R}\frac{1}{\xi-x}\beta^\prime(\xi){\rm d}\xi+\iu \pi \beta^\prime(x), \quad x\in\Omega,
\end{equation}
according to Theorem \ref{theo:extensionofh}.
The function $h$ can then be represented by $\Im\{h\}=\pi\beta^\prime$ as
\begin{equation}\label{eq:Imhtexp}
\Im\{h(x)\}=\sum_{n=1}^{N}c_n\left[p_n(x) +p_n(-x) \right],
\end{equation}
and
\begin{equation}\label{eq:Rehtexp}
\Re\{h(x)\}=b_1x+\sum_{n=1}^{N}c_n\left[\hat{p}_n(x) -\hat{p}_n(-x) \right],
\end{equation}
where $x\in\Omega\subset\R$, $p_n(x)$ are the B-spline basis functions of fixed polynomial order $m-1$,
$b_1\geq 0$ and $c_n\geq 0$ the corresponding optimization variables for $n=1,\ldots,N$, $N$ the number of B-splines
and where $\hat{p}_n(x)$ is the (negative) Hilbert transform \cite{King2009} of the B-spline functions.
Note that the expansions in \eref{eq:Imhtexp} and \eref{eq:Rehtexp} are expressed under the assumption that the imaginary part $\Im\{h(x)\}$ is a symmetric function.
 
As an example, a 2nd order (piecewise linear, ``roof-top'') B-spline 
on a uniform partition is given by $p_n(x)=p(x-n\Delta)$, where
\begin{equation}\label{eq:Approxtriangularpulse}
  p(x)=\left\{\begin{array}{ll}
    1-|x|/\Delta & |x|\leq \Delta, \\ 
    0 & |x|>\Delta, 
  \end{array} \right.
\end{equation}
for $x\in\R$, $\Delta>0$ and its (negative) Hilbert transform $\hat{p}_n(x)=\hat{p}(x-n\Delta)$, where
\begin{equation}\label{eq:ApproxtriangularpulseHilbert}
  \hat{p}(x)
  =\frac{1}{\pi \Delta}\left(2x\ln|x|-(x-\Delta)\ln|x-\Delta|-(x+\Delta)\ln|x+\Delta|\right),
\end{equation}
for $x\in\R$, see \eg \cite{King2009}. Note that the logarithmic singularities in \eref{eq:ApproxtriangularpulseHilbert} cancel, and the function $\hat{p}(x)$ is  continuous on $\R$.
Explicit formulas for general non-uniform B-splines and their Hilbert transforms are given in \cite{Ivanenko+etal2017}.

The convex optimization problem \cite{Boyd+Vandenberghe2004} related to \eref{eq:cvxdefsummary} can now be formulated as
\begin{eqnarray}\label{eq:cvxdef}
\begin{array}{llll}
	& \minimize & & \|h - F \|_{\mrm{L}^p(w,\Omega)}  \\    
	& \subto & &  c_n \geq 0, \quad n=1,\ldots,N,\\  
	&       &  & b_1 \geq 0,
\end{array}
\end{eqnarray}
where $F$ is the given target function, $h$ is given by \eref{eq:Imhtexp} and \eref{eq:Rehtexp} 
and where the minimization is over the variables $b_1$ and $c_n$.

It is finally noted that the uniform continuity of all functions involved implies that the solution to
\eref{eq:cvxdef} can be approximated within an arbitrary accuracy by discretizing the approximation domain $\Omega$ (\ie the computation of the norm) using only a finite number of sample points.
The corresponding numerical problem \eref{eq:cvxdef} can now be solved efficiently by using the CVX Matlab software for disciplined convex programming \cite{Grant+Boyd2012}.
The convex optimization formulation \eref{eq:cvxdef} offers a great advantage in the flexibility in which additional or alternative convex constraints
and formulations can be implemented, see also \cite{Nordebo+etal2014b,Nordebo+etal2017a,Ivanenko+Nordebo2016a,Ivanenko+etal2017}. 

\section{Numerical Examples}
In the following numerical examples, a passive approximation problem is studied in relation to the Example \ref{ex:approxmaxnorm}.
The problem is concerned with the passive realization of an electromagnetic constitutive relation with metamaterial characteristics 
specified over a given finite frequency bandwidth, see also \cite{Gustafsson+Sjoberg2010a}.
Here, the corresponding Herglotz function $h(z)$ is regarded as an holomorphic Fourier transform \cite{Rudin1987} with the variable $x=\Re\{z\}$ 
corresponding to the angular frequency, commonly denoted as $\omega$ (in\unit{rad/s}).

Consider a dielectric metamaterial with a constant, real-valued target permittivity $\epsilon_{\rm t}<0$ to be approximated over an interval $\Omega$, 
and where $F(x)=x\epsilon_{\rm t}$, $h_0(z)=-F(z)$ and hence $b_1^0=-\epsilon_{\rm t}>0$, \cf Example \ref{ex:approxmaxnorm}. 
Let $\Omega=\omega_0[1-\frac{B}{2},1+\frac{B}{2}]$, where $\omega_0$ is the center frequency and $B$ the relative bandwidth with $0<B<2$.
The resulting physical bound can be obtained from \eref{eq:nontrivialbound} as follows
\begin{equation}\label{eq:sumruleconstraint}
\|\epsilon-\epsilon_{\rm t}\|_\infty=\|h-F\|_{\mrm{L}^\infty(w,\Omega)}\geq \frac{\|h-F\|_{\infty}}{\omega_0(1+\frac{B}{2})}
\geq \frac{(\epsilon_{\infty}-\epsilon_{\rm t})\frac{1}{2}B}{1+\frac{B}{2}},
\end{equation}
where $\epsilon$ is the permittivity function of the approximating passive dielectric material, $h(z)=z\epsilon(z)$ the corresponding Herglotz function
and $b_1=\epsilon_{\infty}$ the assumed high-frequency permittivity of the metamaterial, \cf \cite{Gustafsson+Sjoberg2010a}.
Note that here $\|\epsilon-\epsilon_{\rm t}\|_\infty=\|h-F\|_{\mrm{L}^\infty(w,\Omega)}$, where the weight function is $w(x)=\frac{1}{x}$ for $x\in\Omega$, and where $0\notin\Omega$.

Following the derivation given in \cite{Gustafsson+Sjoberg2010a}, it is noted
that the resulting bound given by \eref{eq:sumruleconstraint} is not valid in the case when
the target permittivity $\epsilon_{\rm t}$ has an imaginary part so that $\epsilon_{\rm t}=\Re\{\epsilon_{\rm t}\}+\iu \Im\{\epsilon_{\rm t}\}$ with $\Im\{\epsilon_{\rm t}\}>0$.
However, a straightforward application of the triangle inequality 
\begin{equation}
\left|\epsilon-\Re\{\epsilon_{\rm t}\} \right|=\left|\epsilon-\epsilon_{\rm t}+\iu\Im\{\epsilon_{\rm t}\} \right|
\leq \left|\epsilon-\epsilon_{\rm t} \right|+\Im\{\epsilon_{\rm t}\}, 
\end{equation}
yields the following (Sum rule) bound that is useful when $\Im\{\epsilon_{\rm t}\}$ is small
\begin{equation}\label{eq:sumruleconstraintcomplex}
\|\epsilon-\epsilon_{\rm t}\|_\infty\geq 
\max\left\{\frac{(\epsilon_{\infty}-\Re\{\epsilon_{\rm t}\})\frac{1}{2}B}{1+\frac{B}{2}}- \Im\{\epsilon_{\rm t}\},0 \right\}.
\end{equation}

In the following numerical examples, a dielectric metamaterial is considered with a constant, 
complex valued target permittivity $\epsilon_{\rm t}$ where $\Re\{\epsilon_{\rm t}\}<0$ and $\Im\{\epsilon_{\rm t}\}>0$,
and which is constrained to have a prescribed instantaneous response $\epsilon_\infty=1$.
A physical application for this formulation is with the optimal plasmonic resonances in lossy surrounding media defined in \cite{Nordebo+etal2017a}.
The convex optimization formulation for this problem is 
\begin{equation}\label{eq:cvxdefcomplexmetamaterial}
\begin{array}{llll}
	& \minimize & & \|h-F\|_{\mrm{L}^\infty(w,\Omega)}  \\    
	& \subto & & c_n \geq 0, \quad n=1,\ldots,N, \\  
	&       &  & a_{-1} \leq 0,  
\end{array}
\end{equation}
where $c_n$ and $a_{-1}$ are the optimization variables and the weight function $w$ is defined as above.
Here, the approximating Herglotz function $h(x)=x\epsilon(x)$ is given by
\begin{equation}\label{eq:herglotzregular1}
h(x) =  \frac{a_{-1}}{x}+h_1(x),
\end{equation}
where $h_1(x)$ is represented as in \eref{eq:Imhtexp} and \eref{eq:Rehtexp} with $b_1=\epsilon_\infty=1$, and where $-a_{-1}$ is the 
amplitude of an assumed point mass at $x=0$. In the numerical implementation, $\omega_0=1$ and $N$ uniform linear B-splines are used over the bandwidth $[1-2B,1+2B]$.
The norm is evaluated based on a uniform sampling of $\Omega$ with a fixed sampling  density given by $\Delta x=4B/1000$.
Hence, there are $1000$ sampling points in $[1-2B,1+2B]$ and $250$ sampling points in $\Omega=[1-B/2,1+B/2]$ as $N$ varies between $N=20$ and $N=500$.

\begin{figure}[htb]
\begin{center}
\includegraphics[width=0.95\textwidth]{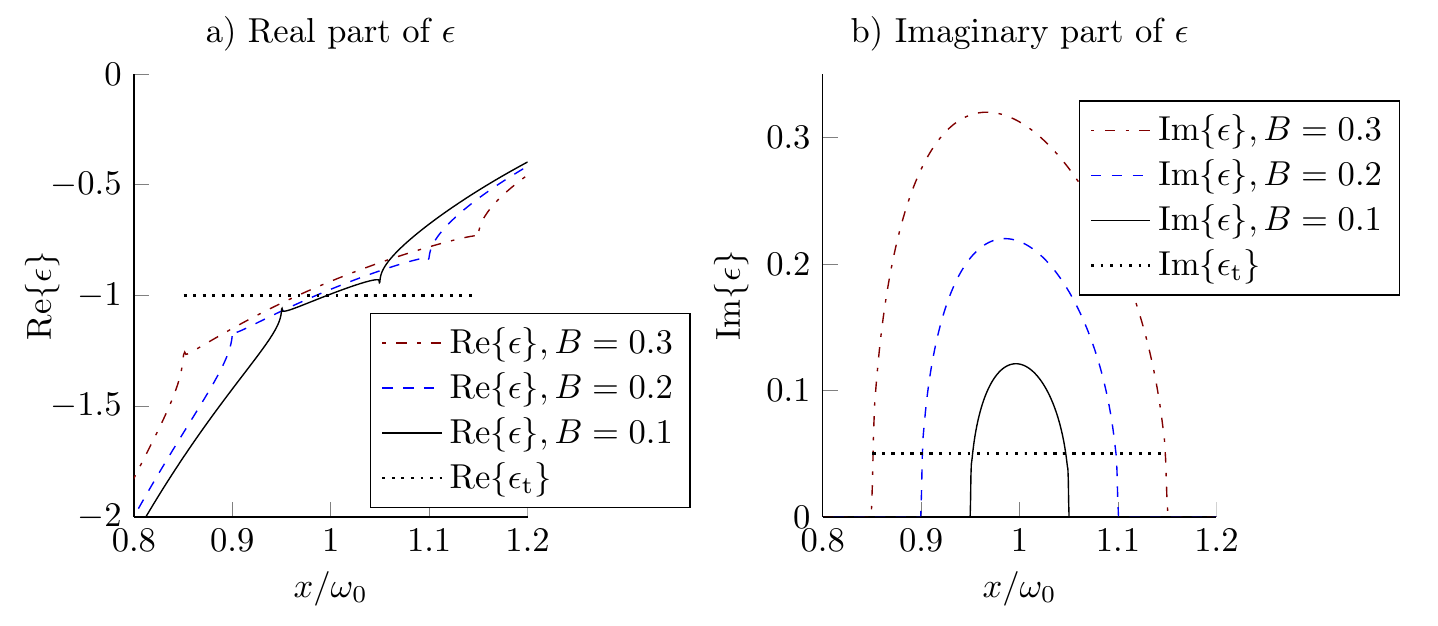}
\end{center}
\vspace{-5mm}
\caption{Optimized passive realization of a lossy metamaterial where the target permittivity is $\epsilon_{\rm t}=-1+\iu0.05$, 
and the approximation interval $\Omega=\omega_0[1-\frac{B}{2},1+\frac{B}{2}]$.}
\label{fig:matfig103-104}
\end{figure}

In Figures \ref{fig:matfig103-104}a-b are illustrated an optimized passive realization of $\epsilon$ using $N=500$ uniform linear B-splines, 
and where the lossy metamaterial is represented by the target permittivity $\epsilon_{\rm t}=-1+\iu0.05$.  
Here, the relative frequency bandwidths are $B=0.1$, $B=0.2$ and $B=0.3$ and the corresponding results are plotted as solid, dashed and dash-dotted lines, respectively, 
and the dotted lines correspond to the target permittvity $\epsilon_{\rm t}$. As seen in the  Figure \ref{fig:matfig103-104}b, the optimization yielded a measure that was supported only
in $\Omega$, except for the point mass at $x=0$ (not shown in the plot).

\begin{figure}[ht!]
\begin{center}
\includegraphics[width=0.95\textwidth]{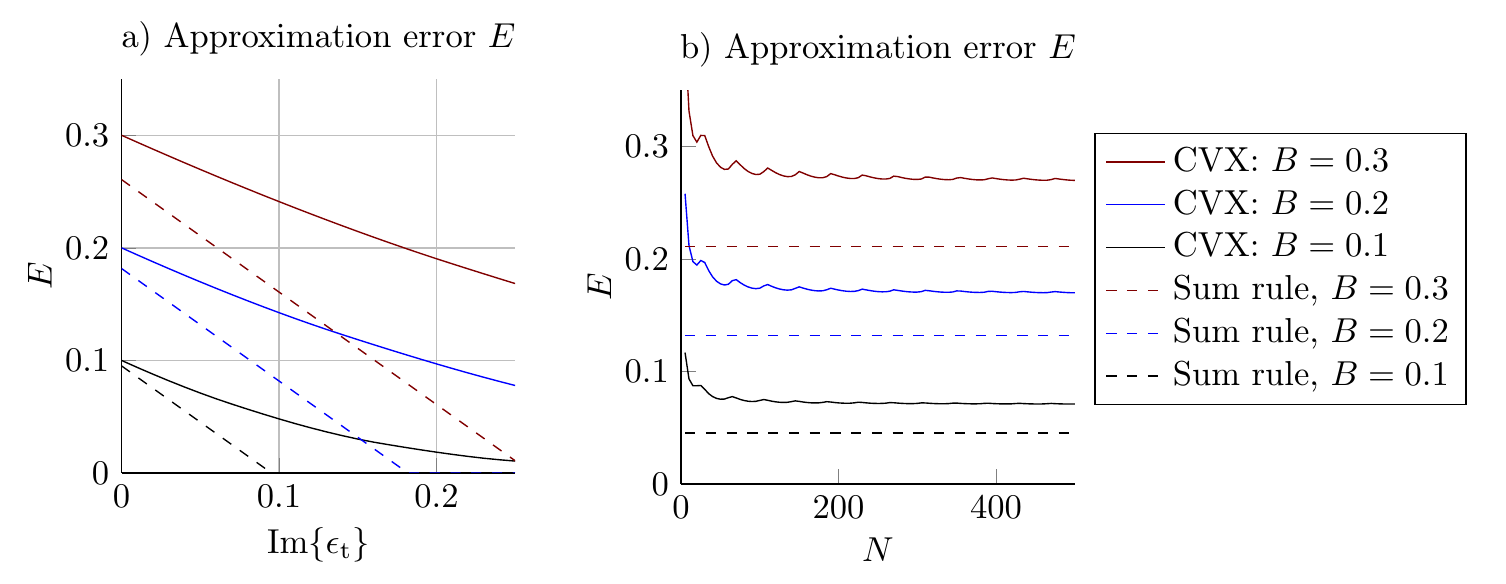}
\end{center}
\vspace{-5mm}
\caption{Approximation error $E$ plotted as a function of a) the imaginary part of the target dielectric permittivity $\Im\{\epsilon_{\rm t}\}$, where $\Re\{\epsilon_{\rm t}\}=-1$ is fixed and $N=500$,
and b) the number of basis functions $N$, where $\epsilon_{\rm t}=-1+\iu0.05$.}
\label{fig:matfig106}
\end{figure}

In Figure \ref{fig:matfig106}a is illustrated the relation between the approximation error $E=\| \epsilon-\epsilon_{\rm t} \|_\infty$ and the imaginary part of the
target permittivity $\Im\{\epsilon_{\rm t}\}$, and where $\Re\{\epsilon_{\rm t}\}=-1$ is fixed. 
In Figure \ref{fig:matfig106}b is illustrated the relation between the approximation error $E$ and the number of basis functions $N$ in the case where $\epsilon_{\rm t}=-1+\iu0.05$. 
The solid and dashed lines correspond to the approximation results (CVX) of \eref{eq:cvxdefcomplexmetamaterial} and the (Sum rule) lower bounds \eref{eq:sumruleconstraintcomplex}, respectively. 
Note that the approximation error increases with increasing relative bandwidth $B$. 
It is seen that the CVX solutions converge well for $N>100$, but do not meet the corresponding sum rule lower bounds as $N$ increases.
This is just a reflection of the fact that the sum rule bounds in \eref{eq:sumruleconstraint} and \eref{eq:sumruleconstraintcomplex} are not tight.
However, as the figure illustrates, the sum rule lower bounds do become tighter as the bandwidth $B$ decreases.
Note finally that according to the Theorem \ref{theo:cvxoptsummary}, the approximation results of \eref{eq:cvxdefcomplexmetamaterial} are asymptotically tight as $N\rightarrow\infty$.

\section{Summary and conclusions}
In this paper, a mathematical framework for passive approximation has been presented.
The objective is to approximate an arbitrary continuous function defined on a finite union of closed and bounded intervals of the real axis (the approximation domain), 
based on the set of Herglotz functions that have H\"{o}lder continuous density in an arbitrary neighborhood of the approximation domain.
The norm used is induced by any weighted $\mrm{L}^p$-norm where $1\leq p\leq\infty$. 
It is proved that the greatest lower bound on the norm of the approximation error can be approached within arbitrary accuracy by using 
a finite uniform B-spline expansion of arbitrary order (linear, quadratic, cubic, etc.) as a generating measure for the approximating Herglotz function.
A numerical example concerning the passive approximation of a metamaterial is given to illustrate the theory. 

\ack
This work was supported by the \href{http://www.stratresearch.se/}{Swedish Foundation for Strategic Research (SSF)} under the program Applied Mathematics and the project \href{http://www.eit.lth.se/index.php?puid=175&projectpage=projektfakta}{Complex analysis and convex optimization for EM design}.

\appendix

\section{Derivation of a non-trivial lower bound}\label{sect:nontrivialbound}

To prove \eref{eq:nontrivialbound}, the technique described in \cite{Gustafsson+Sjoberg2010a} is followed 
based on a composition using the auxiliary Herglotz function $h_\Delta(z)$, defined by the square pulse with
$p_{\Delta}(x)=1$ for $|x|\leq\Delta$, and its associated sum rules. 
The auxiliary Herglotz function and its asymptotics are given by
\begin{equation}
h_\Delta(z)=\frac{1}{\pi}\int_{-\Delta}^{\Delta}\frac{1}{\xi-z}{\rm d}\xi=\frac{1}{\pi}\ln\frac{z-\Delta}{z+\Delta}
=\left\{\begin{array}{l}
 \iu + o(1) \quad \textrm{as}\ z\toh 0, \vspace{0.2cm}\\
 \displaystyle\frac{-2\Delta}{\pi z}+o(z^{-1}) \quad \textrm{as}\ z\toh \infty,
\end{array}\right.
\end{equation}
where $\Im\{h_\Delta(z)\}\geq \frac{1}{2}$ for $|z|\leq \Delta$ and $\Im\{z\}\geq 0$.

Next, the composite Herglotz function ${h}_1(z)$ is defined by 
\begin{equation}
{h}_1(z)=h_\Delta(h(z)+h_0(z)),
\end{equation}
where $h(z)+h_0(z)= (b_1+b_1^0)z+o(z)$ as $z\toh \infty$, yielding the asymptotics
\begin{equation}
{h}_1(z)=
\left\{\begin{array}{l}
 o(z^{-1}) \quad \textrm{as}\ z\toh 0, \vspace{0.2cm}\\
 \displaystyle \frac{-2\Delta }{\pi(b_1+b_1^0)}z^{-1}
 +o(z^{-1}) \quad \textrm{as}\ z\toh \infty.
\end{array}\right.
\end{equation}
The sum rule \eref{eq:Herglotzidentity} for $k=0$ is given by
\begin{equation}
\frac{2}{\pi}\int_{0+}^{\infty}\Im\{h_1(\xi)\}{\rm d} \xi=  a_{-1}-b_{-1}=\frac{2\Delta }{\pi(b_1+b_1^0)}.
\end{equation}
Let $\Delta=\sup_{x\in\Omega}|h(x)+h_0(x)|$, then the following integral inequalities follow
\begin{equation}
\hspace{-1.5cm} \frac{1}{\pi}|\Omega|\leq\frac{2}{\pi}\int_{\Omega}\underbrace{\Im\{h_1(\xi)\}}_{\geq\frac{1}{2}}
{\rm d}\xi\leq\frac{2}{\pi}\int_{0+}^{\infty}\Im\{h_1(\xi)\}{\rm d} \xi
=\frac{2\sup_{x\in\Omega}|h(x)+h_0(x)|}{\pi(b_1+b_1^0)},
\end{equation}
or
\begin{equation}\label{eq:nontrivialbound2}
\|h+h_0\|_\infty\geq (b_1+b_1^0)\frac{1}{2}|\Omega|,
\end{equation}
where $|\Omega|=\int_{\Omega}{\rm d}\xi$.


\end{document}